\documentclass{article}
\usepackage[latin1]{inputenc}
\usepackage{amsmath}
\usepackage{amsthm}
\usepackage{amsfonts} 

\makeatletter
\def\multilimits@{\bgroup
  \Let@
  \restore@math@cr
  \default@tag
 \baselineskip\fontdimen10 \scriptfont\tw@
 \advance\baselineskip\fontdimen12 \scriptfont\tw@
 \lineskip\thr@@\fontdimen8 \scriptfont\thr@@
 \lineskiplimit\lineskip
 \vbox\bgroup\ialign\bgroup\hfil$\m@th\scriptstyle{##}$\hfil\crcr}
\def\Sb{_\multilimits@}
\def\Sp{^\multilimits@}
\def\endSb{\crcr\egroup\egroup\egroup}

\makeatother

\newif\iffrac
\fractrue

\theoremstyle{plain}
\newtheorem{lemma}{Lemma}
\newtheorem{proposition}[lemma]{Proposition}
\newtheorem{theorem}{Theorem}
\newtheorem{corollary}[lemma]{Corollary}

\theoremstyle{remark}
\newtheorem*{remark}{Remark}

\newlength{\fixboxwidth}
\renewcommand{\phi}{\varphi}
\renewcommand{\epsilon}{\varepsilon}
\renewcommand{\rho}{\varrho}

\begin{document}

\title{The average distance property of\\classical Banach spaces II}
\author{Aicke Hinrichs\thanks{supported by DFG grant Hi 584/2-2} \ and
J{\"o}rg Wenzel\thanks{supported by DFG grant We 1868/1-1
 }}

\maketitle

\begin{abstract}
  A Banach space X has the average distance property (ADP) if
  there exists a unique real number $r$ such that for each
  positive integer $n$ and all $x_1,\ldots,x_n$ in the unit sphere
  of $X$ there is some $x$ in the unit sphere of $X$ such that
  \[ \frac{1}{n} \sum_{k=1}^n \|x_k-x\|=r.\]
  We show that $l_p$  does not have the average distance
  property if $p>2$. This completes the study of the ADP for
  $l_p$ spaces.
\end{abstract}


\section{Introduction}
\label{sec:intro}

The aim of this note is to finish the study of the average
distance property of $l_p$ and $L_p[0,1]$ for $1\le p\le \infty$
using and refining the method introduced in \cite{hinrichs00:_l}.
We start giving a short review of that method. The reader is
referred to \cite{hinrichs00:_l} for further information and to
the pointers to the literature therein.

A rendezvous number of a metric space $(M,d)$ is a real number $r$
with the property that for each positive integer $n$ and
$x_1,\ldots,x_n\in M$ there exists $x\in M$ such that
\[ \frac{1}{n}\sum_{k=1}^n d(x_k,x) = r.\]
We say that a (real or complex) Banach space $X$ has the average
distance property (ADP for short) if its unit sphere has a unique
rendezvous number.
It is known that $l_2$ and $L_2[0,1]$ have the ADP
\cite{morris83:_average} and that $l_p$ and $L_p[0,1]$ do not have the
ADP if $1\le p<2$ and if $p\ge 3$, see \cite{lin97:_average} and
\cite{hinrichs00:_l}, respectively. Here we prove the following
result.
\begin{theorem}
  \label{thm:main} For $p>2$, $l_p$ and $L_p[0,1]$ do not have the
  ADP.
\end{theorem}

In \cite{hinrichs00:_l}, using an improved Clarkson inequality, the
study of the ADP for $l_p$ and $L_p$ in the case $p>2$ was reduced to
the study of a scalar function as follows. For $n\in \mathbb N$, $p>2$
and $x,y_1,\dots,y_n\in l_p$ or $L_p$ such that $\|x\|^p=1/n$ and
$\sum_{i=1}^n\|y_i\|^p=1$ define
\begin{equation}\label{eq:5}
  \sigma_i:=\frac{\|x-y_i\|^p}{\|x\|^p+\|y_i\|^p}
  \quad\mbox{and}\quad
  \alpha_i:=\frac{\|x\|^p+\|y_i\|^p}2.
\end{equation}
It follows that
\[ \frac1{2n}\leq \alpha_i\leq \frac{n+1}{2n}
\quad\mbox{and}\quad
\sum_{i=1}^n\alpha_i=1.
\]
Let $u_i\in[-1,+1]$ be defined by the relation
\begin{equation}\label{eq:4}
  \sigma_i=\frac{(1-u_i)^p}{1+|u_i|^p}.
\end{equation}
and let
\[ \phi(x,y_1,\dots,y_n) := 2^{-n}
\sum_{\epsilon_1,\dots,\epsilon_n=\pm1} \Big(
\sum_{i=1}^n \alpha_i \frac{(1+\epsilon_i
  u_i)^p}{1+|u_i|^p}\Big)^{1/p}.
\]
As pointed out in \cite{hinrichs00:_l}, in order to prove
Theorem~\ref{thm:main} for a fixed $p>2$, it suffices to find $n$ such
that $\phi>1$ for $(u_1,\ldots,u_n)\not=(0,\ldots,0)$.

Considering the case $u_i=1$ for $i=1,\dots,n$, $\alpha_i=1/(2n)$ for
$i=1,\ldots,n-1$, and $\alpha_n=(n+1)/(2n)$ yields that
\begin{eqnarray*}
  \lefteqn{
    2^{-n} \sum_{\epsilon_1,\ldots,\epsilon_n=\pm1}\Big(
    \sum_{i=1}^n \alpha_i \frac{(1+\epsilon_i
      u_i)^p}{1+|u_i|^p}\Big)^{1/p}
    } \hspace{1cm} \\
  & = &
  2^{1-2/p} 2^{-n} \sum_{\epsilon_1,\ldots,\epsilon_n=\pm1}\Big(
  1+\frac1{2n} \sum_{i=1}^{n-1}\epsilon_i + \frac{n+1}{2n}\epsilon_n
  \Big)^{1/p} \\
  & = & 2^{-2/p} \sum_{\epsilon_n=\pm1} \Bigg( 2^{-n+1}
  \sum_{\epsilon_1,\ldots,\epsilon_{n-1}=\pm1} \Big(
  1+\frac1{2n}\sum_{i=1}^{n-1} \epsilon_i + \frac{n+1}{2n}\epsilon_n
  \Big)^{1/p} \Bigg) \\
  &\leq&
  2^{-2/p} \sum_{\epsilon_n=\pm1}
  \Big(
  1+ 2^{-n+1} \sum_{\epsilon_1,\ldots,\epsilon_{n-1}=\pm1}
  \frac1{2n}\sum_{i=1}^{n-1} \epsilon_i + \frac{n+1}{2n} \epsilon_n
  \Big)^{1/p}
  \\
  & = &
  2^{-2/p} \Bigg( \Big( \frac32+\frac1n\Big)^{1/p} +
  \Big(\frac12-\frac1n\Big)^{1/p} \Bigg) \\
  &\leq&
  2^{-2/p} \Bigg( \Big( \frac32\Big)^{1/p} +
  \Big(\frac12\Big)^{1/p} \Bigg)
  =
  \frac{3^{1/p}+1}{8^{1/p}}
\end{eqnarray*}
which is smaller than $1$ for $p<2.10528\dots$

This shows that, in contrast to \cite{hinrichs00:_l}, we have to take
into account the concrete definition of the $u_i$'s and $\alpha_i$'s
to be able to cover also the cases where $p$ is close to $2$. This
will be done in Proposition~\ref{lem:4}.

The remaining part of the paper is the proof of
Theorem~\ref{thm:main}, which follows from the upcoming
Propositions \ref{prop:1} and \ref{prop:2}.


\section{The relation of $\alpha_i$ and $u_i$}
\label{sec:relation}

We begin by providing an auxiliary estimate.

\begin{lemma}\label{lem:5}
  \[ \frac{1+u}{(1+u^p)^{1/p}} \geq 1+(2^{1-1/p}-1)u
  \]
  for $u\in[0,1]$.
\end{lemma}
\begin{proof}
  Let $g(u):=(1+u)/(1+u^p)^{1/p}$. Note that
  \[ g' (u)= \frac{1-u^{p-1}}{(1+u^p)^{1+1/p}} \geq 0
  \]
  while
  \[ g''(u)=-\frac{(p+1)(1-u^{p-1})u^{p-1}+(p-1)(1+u^p)u^{p-2}}
  {(1+u^p)^{2+1/p}} \leq 0.
  \]
  This means that $g$ is a concave function on $[0,1]$ and therefore
  $g(u)\geq g(0)+(g(1)-g(0))u$. This proves the assertion.
\end{proof}

\begin{proposition}\label{lem:4}
  If $\alpha_i$ and $u_i$ are defined by~{\rm (\ref{eq:5})}
  and~{\rm (\ref{eq:4})}, then
  \[
  |u_i| \leq c_1 n^{-1/p}\alpha_i^{-1/p},
  \]
where $c_1=\max(2^{1-1/p},1/(2-2^{1/p}))$.
\end{proposition}
\begin{proof}
We split the proof into three cases.

  \textbf{First case:}
  \[ \alpha_i \leq \frac1n.
  \]
  Since $c_1\geq1$, in this case
  \[ |u_i| \leq 1 \leq
  (n\alpha_i)^{-1/p}\leq c_1 n^{-1/p}\alpha_i^{-1/p}.
  \]

  \textbf{Second case:}
  \[ \alpha_i>\frac1n
  \quad\mbox{and}\quad
  u_i\geq0.
  \]
  Then
  \[ \frac{\|y_i\|}{(2\alpha_i)^{1/p}} =
  \Big(1-\frac1{2\alpha_in}\Big)^{1/p} >
  \Big(\frac1{2\alpha_in}\Big)^{1/p} = \frac{\|x\|}{(2\alpha_i)^{1/p}}
  \]
  and it follows from the definition~(\ref{eq:4}) of $u_i$ that
  \begin{eqnarray*}
    (1-u_i)^p \geq \frac{(1-u_i)^p}{1+u_i^p} = \sigma_i & = &
    \frac{\|x-y_i\|^p}{\|x\|^p+\|y_i\|^p}  =
    \Big\|\frac x{(2\alpha_i)^{1/p}}
    -\frac{y_i}{(2\alpha_i)^{1/p}}\Big\|^p\\
    & \geq &
    \Bigg(\Big(1-\frac1{2\alpha_in}\Big)^{1/p} -
    \Big(\frac1{2\alpha_in}\Big)^{1/p}\Bigg)^p.
  \end{eqnarray*}
  Now, using the relations
  \[ 1-\frac1{2\alpha_in} \leq
  \Big(1-\frac1{2\alpha_in}\Big)^{1/p}
  \quad\mbox{and}\quad
  \frac1{2\alpha_in} \leq
  \Big(\frac1{2\alpha_in}\Big)^{1/p}
  \]
  which follow from $\alpha_i\geq1/(2n)$, we obtain
  \[ u_i \leq 1- \Big(1-\frac1{2\alpha_in}\Big)^{1/p} +
  \Big(\frac1{2\alpha_in}\Big)^{1/p} \leq
  \frac1{2\alpha_in} +
  \Big(\frac1{2\alpha_in}\Big)^{1/p} \leq
  2\Big(\frac1{2\alpha_in}\Big)^{1/p}.
  \]
  From this we get
  \[ |u_i| = u_i \leq
  2^{1-1/p}n^{-1/p}\alpha_i^{-1/p}.
  \]
  \textbf{Third case:}
  \[ \alpha_i>\frac1n
  \quad\mbox{and}\quad
  u_i<0.
  \]
  It follows from Lemma~\ref{lem:5} for $u=-u_i$ that
  \[ 1-(2^{1-1/p}-1)u_i \leq \sigma_i^{1/p} \leq
  \Big(1-\frac1{2\alpha_in}\Big)^{1/p} +
  \Big(\frac1{2\alpha_in}\Big)^{1/p} \leq
  1+\Big(\frac1{2\alpha_in}\Big)^{1/p}.
  \]
  Finally in this case
  \[ |u_i|= -u_i\leq \frac1{2-2^{1/p}} n^{-1/p} \alpha_i^{-1/p}.
  \]
\end{proof}

With this proposition in hand, we can forget about the concrete nature
of the $\alpha_i$'s and $u_i$'s. All we have to show is that for given
$n$ and $\alpha_1,\dots,\alpha_n$ such that
\[ \frac1{2n}\leq\alpha_i\leq\frac{n+1}{2n}
\quad\mbox{and}\quad
\sum_{i=1}^n\alpha_i=1
\]
the function
\[ \phi(u_1,\ldots,u_n):=2^{-n}\sum_{\epsilon_1,\ldots,\epsilon_n=\pm1}
\Big(\sum_{i=1}^n \alpha_i
\frac{(1+\epsilon_iu_i)^p}{1+|u_i|^p}\Big)^{1/p}
\]
is bigger than one as long as
\begin{equation}\label{eq:7}
  |u_i|\leq c_1n^{-1/p}\alpha_i^{-1/p}
\end{equation}
and $(u_1,\ldots,u_n)\not=(0,\ldots,0)$.

Since all relations on the $u_i$'s are symmetric and since the
function $\phi$ is symmetric in $u_i$, we can henceforth assume that
$u_i\geq0$.


\section{Proof of $\phi>1$, the case of many large $u_i$'s}
\label{sec:proof}

\begin{corollary}\label{cor:1}
  \[ \Big(\sum_{i=1}^n (\alpha_i u_i)^2\Big)^{1/2} \leq c_1n^{-1/p}.
  \]
\end{corollary}
\begin{proof}
  It follows from~(\ref{eq:7}) that
  \[ \Big(\sum_{i=1}^n (\alpha_iu_i)^2\Big)^{1/2} \leq
  c_1 n^{-1/p} \Big(\sum_{i=1}^n \alpha_i^{2-2/p}\Big)^{1/2}.
  \]
  Since $2-2/p>1$ and $\alpha_i<1$ we have
  \[ \sum_{i=1}^n\alpha_i^{2-2/p} \leq \sum_{i=1}^n\alpha_i = 1,
  \]
  which proves the assertion.
\end{proof}

\begin{lemma}
  \label{lem:2}
  We have
  \begin{equation}\label{eq:2}
    v(u):=\frac{(1+u)^p+(1-u)^p}{2(1+u^p)}\geq 1+c_2 u^p
  \end{equation}
  and
  \begin{equation}\label{eq:1}
    w(u):=\frac{(1+u)^p-(1-u)^p}{2(1+u^p)}\leq c_3u
  \end{equation}
  for $u\in[0,1]$, where $c_2:=2^{p-2}-1$ and $c_3:=p2^{p-1}$.
\end{lemma}
\begin{proof}
  To see~(\ref{eq:2}), we let
  \[ g(u):=\frac{(1+u)^p+(1-u)^p-2}{u^p}
  \]
  and use the fact that $(1+u)^{p-1}+(1-u)^{p-1}$ is non-increasing
  for $p>2$, to compute
  \[
  g'(u)  =  \frac p{u^{p+1}} (2-(1+u)^{p-1}-(1-u)^{p-1}) \leq 0.
  \]
  Therefore $g(u)\geq g(1) = 2^p-2$, which yields
  \[ (1+u)^p+(1-u)^p\geq 2+(2^p-2)u^p = 2(1+u^p)+(2^p-4)u^p.
  \]
  Division by $2(1+u^p)$ and $1+u^p\leq2$ proves~(\ref{eq:2}).

  Since $2u/(1+u)\leq 1$, Bernoulli's inequality states
  \[ \frac{(1-u)^p}{(1+u)^p} =
  \Big(1-\frac{2u}{1+u}\Big)^p \geq 1-\frac{2pu}{1+u}.
  \]
  It follows that
  \[ \frac{(1+u)^p-(1-u)^p}{1+u^p}
  = \frac{(1+u)^p}{1+u^p}\Big(1-\frac{(1-u)^p}{(1+u)^p}\Big)
  \leq 2pu\frac{(1+u)^{p-1}}{1+u^p} \leq p2^pu
  \]
  which proves~(\ref{eq:1}).
\end{proof}

The following Lemma is known as a subgaussian tail estimate for
Rademacher averages and is by now classical. A proof can be found
e.~g.~in \cite[p.~90]{ledoux91:_probab_banac}.
\begin{lemma}\label{lem:3}
  For a given vector $x=(\xi_1,\ldots,\xi_n)$, let
  $\|x\|_2:=\Big(\sum_{i=1}^n |\xi_i|^2\Big)^{1/2}$ and $\mathbb
  B:=\{(\epsilon_1,\ldots,\epsilon_n) :
  \sum_{i=1}^n\epsilon_i\xi_i>t\|x\|_2\}$, then
  \[
  2^{-n} |\mathbb B| \leq
  e^{-t^2/2}.
  \]
\end{lemma}

We are now ready to tackle the case, where `many' of the $u_i$'s are
bigger than~$1/2$.

\begin{proposition}
\label{prop:1}
  There exists $n_1$ such that for all $n>n_1$ we have
  \[ \phi(u_1,\ldots,u_n)>1
  \]
  if $|\mathbb{A}|>n/2$, where $\mathbb{A}:=\{i:u_i>1/2\}$.
\end{proposition}
\begin{proof}
  With $v$ and $w$ defined as in Lemma~\ref{lem:2}, observe that
  \[ v(u)+\epsilon w(u) = \frac{(1+\epsilon u)^p}{1+u^p}
  \]
  for $\epsilon=\pm1$.  Put
  \[ \mathbb B:=\Big\{(\epsilon_1,\ldots,\epsilon_n):
  -\sum_{i=1}^n\alpha_i \epsilon_i w(u_i) \leq
  (2\log n)^{1/2} c_3c_1n^{-1/p}
  \Big\}.
  \]
  Since by~(\ref{eq:1}) and Corollary~\ref{cor:1}
  \[ \Big(\sum_{i=1}^n \big( \alpha_i
  w(u_i) \big)^2\Big)^{1/2}
  \leq c_3 \Big(\sum_{i=1}^n \big(\alpha_iu_i\big)^2\Big)^{1/2}
  \leq  c_3c_1n^{-1/p}
  \]
  it follows from Lemma~\ref{lem:3} that
  \[
  2^{-n} |\mathbb B| \geq
  1-\frac1n.
  \]
  With these preliminaries we can estimate $\phi$ as follows
  \begin{eqnarray*}
    \phi(u_1,\ldots,u_n)&\geq&
    2^{-n}\sum_{(\epsilon_1,\ldots,\epsilon_n)\in\mathbb B}
    \Big(\sum_{i=1}^n\alpha_i
    v(u_i) + \sum_{i=1}^n\alpha_i
    \epsilon_i w(u_i) \Big)^{1/p} \\
    & \geq &
    \Big(1-\frac1n\Big)
    \Big(
    \sum_{i=1}^n \alpha_i v(u_i) -
    (2\log n)^{1/2}c_3c_1n^{-1/p}
    \Big)^{1/p}.
  \end{eqnarray*}
  From ~(\ref{eq:2}) and the assumption on $\mathbb{A}$ it follows
  that
  \[ \sum_{i=1}^n \alpha_i v(u_i) \geq
  \sum_{i=1}^n\alpha_i+\sum_{i\in \mathbb{A}} \alpha_ic_2 u_i^p
  \geq 1 + \frac n2 \frac1{2n} c_2 2^{-p} = 1 + c_4,
  \]
  where $c_4:=c_22^{-p-2}$.

  Since $c_4>0$, we can now choose $n_1$ so that for all $n>n_1$
  \[ (2\log n)^{1/2}c_3c_1n^{-1/p} < \frac{c_4}2
  \quad\mbox{and}\quad
  \Big(1-\frac1n\Big)\Big(1+\frac{c_4}2\Big)^{1/p} >
  \Big(1 +\frac{c_4}4\Big)^{1/p}.
  \]
  By these assumptions on $n$
  \begin{eqnarray*}
    \phi(u_1,\ldots,u_n)
    & \geq &
    \Big(1-\frac1n\Big)
    \Big(
    1+c_4 - (2\log n)^{1/2}c_3c_1n^{-1/p}
    \Big)^{1/p} \\
    & \geq &
    \Big(1-\frac1n\Big)
    \Big(
    1+\frac{c_4}2
    \Big)^{1/p} \\
    & \geq &
    \Big(1+\frac{c_4}4\Big)^{1/p}.
  \end{eqnarray*}
  This proves the assertion.
\end{proof}


\section{Proof of $\phi>1$, the case of few large $u_i$'s}
\label{sec:proof_few}

From now on, we will only deal with the case $|\mathbb{A}|\leq n/2$.
So for the rest of this section, we assume that
\begin{equation}\label{eq:8}
  |\mathbb{A}|\leq \frac n2,
  \quad\mbox{where $\mathbb{A}=\{i:u_i>1/2\}$.}
\end{equation}

\begin{lemma}\label{lem:1}
  Denote
\iffrac
  \begin{equation}\label{eq:9}
    f(u):=\frac{(1-u^2)^p}{1+u^p}
    \frac{(1+u^{p-1})^{\frac p{p-1}} - (1-u^{p-1})^{\frac p{p-1}}}
    {(1+u)^p(1-u^{p-1})^{\frac p{p-1}} - (1-u)^p(1+u^{p-1})^{\frac
        p{p-1}}}.
  \end{equation}
\else
  \begin{equation}\label{eq:9}
    f(u):=\frac{(1-u^2)^p}{1+u^p}
    \frac{(1+u^{p-1})^{p/(p-1)} - (1-u^{p-1})^{p/(p-1)}}
    {(1+u)^p(1-u^{p-1})^{p/(p-1)} - (1-u)^p(1+u^{p-1})^{p/(p-1)}}.
  \end{equation}
\fi
  Then $\lim_{u\to0}f(u)=\lim_{u\to1}f(u)=0$ and $f$ is bounded on
  $[0,1]$.
\end{lemma}
\begin{proof}
  Note that the derivative of the function
  \iffrac
  $(1\pm u)^p(1\mp u^{p-1})^{\frac p{p-1}}$ is
  \else
  $(1\pm u)^p(1\mp u^{p-1})^{p/(p-1)}$ is
  \fi
  \iffrac
  \[ \pm p(1\pm u)^{p-1}(1\mp u^{p-1})^{\frac p{p-1}} \mp
  p(1\pm u)^p(1\mp u^{p-1})^{\frac 1{p-1}}u^{p-2}.
  \]
  \else
  \[ \pm p(1\pm u)^{p-1}(1\mp u^{p-1})^{p/(p-1)} \mp
  p(1\pm u)^p(1\mp u^{p-1})^{1/(p-1)}u^{p-2}.
  \]
  \fi
  Since $p>2$ we therefore have
\iffrac
  \[ \lim_{u\to0} \frac d{du} (1+ u)^p(1- u^{p-1})^{\frac p{p-1}}
  - \lim_{u\to0} \frac d{du} (1- u)^p(1+ u^{p-1})^{\frac p{p-1}}
  = 2 p.
  \]
\else
  \[ \lim_{u\to0} \frac d{du} (1+ u)^p(1- u^{p-1})^{p/(p-1)}
  -\lim_{u\to0} \frac d{du} (1- u)^p(1+ u^{p-1})^{p/(p-1)}
  = 2 p.
  \]
\fi
  By l'Hospital's rule
  \begin{eqnarray*}
    \lim_{u\to0}f(u) & = &
\iffrac
    \lim_{u\to0} \frac{(1+u^{p-1})^{\frac p{p-1}} - (1-u^{p-1})^{\frac
        p{p-1}}}
    {(1+u)^p(1-u^{p-1})^{\frac p{p-1}} -
      (1-u)^p(1+u^{p-1})^{\frac p{p-1}}} \\
    & = &
    \lim_{u\to0}\frac{p(1+u^{p-1})^{\frac 1{p-1}}u^{p-2} +
      p(1-u^{p-1})^{\frac 1{p-1}}u^{p-2}}
    {\frac d{du} (1+ u)^p(1- u^{p-1})^{\frac p{p-1}} - \frac d{du} (1-
      u)^p(1+ u^{p-1})^{\frac p{p-1}}}\\
\else
    \lim_{u\to0} \frac{(1+u^{p-1})^{p/(p-1)} - (1-u^{p-1})^{p/(p-1)}}
    {(1+u)^p(1-u^{p-1})^{p/(p-1)} - (1-u)^p(1+u^{p-1})^{p/(p-1)}} \\
    & = &
    \lim_{u\to0}\frac{p(1+u^{p-1})^{1/(p-1)}u^{p-2} +
      p(1-u^{p-1})^{1/(p-1)}u^{p-2}}
    {\frac d{du} (1+ u)^p(1- u^{p-1})^{p/(p-1)} - \frac d{du} (1-
      u)^p(1+ u^{p-1})^{p/(p-1)}}\\
\fi
    & = & 0.
  \end{eqnarray*}

  On the other hand, again by l'Hospital's rule it follows that
\iffrac
  \[ \lim_{u\to1}\frac{(1-u^{p-1})^{\frac 1{p-1}}}{1-u} =
  \lim_{u\to1} \frac{u^{p-2}}{(1-u^{p-1})^{\frac {p-2}{p-1}}} =
  +\infty.
  \]
\else
  \[ \lim_{u\to1}\frac{(1-u^{p-1})^{1/(p-1)}}{1-u} =
  \lim_{u\to1} \frac{u^{p-2}}{(1-u^{p-1})^{(p-2)/(p-1)}} = +\infty.
  \]
\fi
  Therefore
  \begin{eqnarray*}
    \lim_{u\to1}f(u) & = &
\iffrac
    \lim_{u\to1} \frac{2^{\frac
        1{p-1}}(1-u^2)^p}{(1+u)^p(1-u^{p-1})^{\frac p{p-1}} -
      (1-u)^p(1+u^{p-1})^{\frac p{p-1}}} \\
    & = &
    \lim_{u\to1} \frac{2^{\frac 1{p-1}}}{\displaystyle
      \frac{(1-u^{p-1})^{\frac p{p-1}}}{(1-u)^p} -
      \frac{(1+u^{p-1})^{\frac p{p-1}}}{(1+u)^p}} \\
\else
    \lim_{u\to1} \frac{2^{1/(p-1)}(1-u^2)^p}{(1+u)^p(1-u^{p-1})^{p/(p-1)} -
      (1-u)^p(1+u^{p-1})^{p/(p-1)}} \\
    & = &
    \lim_{u\to1} \frac{2^{1/(p-1)}}{\displaystyle
      \frac{(1-u^{p-1})^{p/(p-1)}}{(1-u)^p} -
      \frac{(1+u^{p-1})^{p/(p-1)}}{(1+u)^p}} \\
\fi
    & = & 0.
  \end{eqnarray*}
  The boundedness of $f$ on $[0,1]$ now follows from its continuity in
  $(0,1)$ and the boundedness of the limits of $f(u)$ for $u\to0$ and
  $u\to1$.
\end{proof}

We can now also treat the remaining case, where only `few' of the
$u_i$'s are bigger than~$1/2$. In this case, the next proposition
shows that $\phi(u_1,\ldots,u_n) > \phi(0,\ldots,0)=1$, provided
that $n$ is big enough. This completes the proof of
Theorem~\ref{thm:main}.

\begin{proposition}
  \label{prop:2}
  There exists $n_2 \geq n_1$ such that for all $n>n_2$ we have
  \[ \frac{\partial\phi}{\partial u_j}(u_1,\dots,u_n) > 0
  \]
  for all $j=1,\ldots,n$ and all $u_1,\dots,u_n$
  satisfying~{\rm (\ref{eq:8})}.
\end{proposition}
\begin{proof}
  Note that
  \[ \frac{\partial\phi}{\partial u_j}(u_1,\ldots,u_n) =
  \frac{\alpha_j}{(1+u_j^p)^2} 2^{-n}
  \sum_{\epsilon_1,\ldots,\epsilon_n=\pm1} \epsilon_j
  \frac{(1+\epsilon_ju_j)^{p-1} (1-\epsilon_ju_j^{p-1})}
  {\displaystyle\Big(\sum_{i=1}^n \alpha_i
    \frac{(1+\epsilon_iu_i)^p}{1+u_i^p}\Big)^{1-1/p}}.
  \]
  We will show that for every
  $\epsilon_1,\ldots,\epsilon_{j-1},\epsilon_{j+1},\ldots\epsilon_n$
  the summand
  \[ \sum_{\epsilon_j=\pm1}
  \epsilon_j
  \frac{(1+\epsilon_ju_j)^{p-1} (1-\epsilon_ju_j^{p-1})}
  {\displaystyle\Big(\sum_{i=1}^n \alpha_i
    \frac{(1+\epsilon_iu_i)^p}{1+u_i^p}\Big)^{1-1/p}}
  \]
  is positive.

  To this end we denote
  \[ a_j(u_1,\ldots,u_n):= \sum\begin{Sb}i=1\\i\not=j\end{Sb}^{n}
  \alpha_i \frac{(1+\epsilon_iu_i)^p}{1+u_i^p}
  \]
  and show that
  \[ \frac{(1+u_j)^{p-1} (1-u_j^{p-1})}
  {\displaystyle\Big(a_j(u_1,\ldots,u_n) + \alpha_j
    \frac{(1+u_j)^p}{1+u_j^p}\Big)^{1-1/p}}
  >
  \frac{(1-u_j)^{p-1} (1+u_j^{p-1})}
  {\displaystyle\Big(a_j(u_1,\ldots,u_n) + \alpha_j
    \frac{(1-u_j)^p}{1+u_j^p}\Big)^{1-1/p}}.
  \]
  Some manipulations show that this is equivalent to
  \[ a_j(u_1,\ldots,u_n) > \alpha_j f(u_j),
  \]
  where $f$ is the function defined in~(\ref{eq:9}) in
  Lemma~\ref{lem:1}.

  Using~(\ref{eq:8}), we see that
  \[ a_j(u_1,\ldots,u_n)
  \geq \sum\begin{Sb}i\not\in\mathbb{A}\\i\not=j\end{Sb}
  \alpha_i \frac{(1-u_i)^p}{1+u_i^p}
  \geq \Big(\frac n2 - 1\Big) \frac1{2n} \frac{2^{-p}}{1+2^{-p}} \geq
  \frac18 \frac1{1+2^p}= c_5,
  \]
  if $n\geq4$ and $c_5:=1/(8+2^{p+3})$.  It is hence enough to
  show that
  \begin{equation}\label{eq:3}
    c_5 > \alpha_j f(u_j).
  \end{equation}

  Since $\lim_{u\to0} f(u)=0$ by Lemma~\ref{lem:1}, we can find
  $\delta>0$ small enough such that
  \[ f(u) < c_5
  \]
  for $u^p<\delta$.
  Since $f$ is also bounded by Lemma~\ref{lem:1}, we can choose
  \[ n \geq n_2 := \max\Big(\frac{c_1^p \|f\|_\infty}{c_5\delta}, n_1,
  4\Big).
  \]

  If $\alpha_j < c_5/\|f\|_\infty$ then obviously~(\ref{eq:3})
  holds.

  If on the other hand $\alpha_j\geq c_5/\|f\|_\infty$ then
  \[ \alpha_j n > \frac{c_5}{\|f\|_\infty}
  \frac{c_1^p\|f\|_\infty}{c_5\delta} = \frac{c_1^p}\delta
  \]
  and by~(\ref{eq:7})
  \[ u_j^p\leq\frac{c_1^p}{n\alpha_j} < \delta.
  \]
  Consequently
  \[ \alpha_jf(u_j)<\alpha_jc_5\leq c_5,
  \]
  since $\alpha_j\leq 1$.

  This proves the assertion.
\end{proof}

\begin{remark}
  Using the methods developed in Sections~\ref{sec:proof}
  and~\ref{sec:proof_few}, it can be shown that without
  Relation~(\ref{eq:7}) one can prove the result of the main theorem
  for all $p>p_0$, where
  \[ p_0 := \inf\{p>2: g\geq 2^{(1+1/p)}\} = 2.2751\ldots
  \]
  and
  \[ g(u):=\Big(1+\frac{(1+u)^p}{1+u^p}\Big)^{1/p} +
  \Big(1+\frac{(1-u)^p}{1+u^p}\Big)^{1/p}.
  \]
\end{remark}


\begin{thebibliography}{1}

\bibitem{hinrichs00:_l}
A.~Hinrichs.
\newblock The average distance property of classical {B}anach spaces.
\newblock {\em Bull. Austral. Math. Soc.}, 62(1):119--134, 2000.

\bibitem{ledoux91:_probab_banac}
M.~Ledoux and M.~Talagrand.
\newblock {\em Probability in {Banach} spaces}.
\newblock Springer, Berlin--Heidelberg, 1991.

\bibitem{lin97:_average}
P.-K. Lin.
\newblock The average distance property of {B}anach spaces.
\newblock {\em Arch. Math. (Basel)}, 68(6):496--502, 1997.

\bibitem{morris83:_average}
S.~A. Morris and P.~Nickolas.
\newblock On the average distance property of compact connected metric spaces.
\newblock {\em Arch. Math. (Basel)}, 40(5):459--463, 1983.

\end{thebibliography}

\bibliographystyle{abbrv}

\noindent address: Mathematisches Institut, FSU Jena, D-07743 {\sc
Jena}, Germany\\
e-mail: {\tt nah@rz.uni-jena.de, wenzel@minet.uni-jena.de}

\end{document}